# Fitted Reproducing Kernel Method for Solving a Class of Third-Order Periodic Boundary Value Problems

[1]Asad Freihat, [2]Radwan Abu-Gdairi, [3]Hammad Khalil, [4]Eman Abuteen, [1]Mohammed Al-Smadi and [3]Rahmat Ali Khan

[1]Department of Applied Science, Ajloun College, Al Balqa Applied University, Ajloun 26816, Jordan
[2]Department of Mathematics, Faculty of Science, Zarqa University, Zarqa 13110, Jordan
[3]Department of Mathematics, University of Malakand, Chakadara Dir(L), Khyber Pakhtunkhwa, Pakistan
[4]Department of Applied Science, Faculty of Engineering Technology, Al-Balqa Applied University, Jordan

**Abstract:** In this article, the reproducing kernel Hilbert space $W_2^4[0, 1]$ is employed for solving a class of third-order periodic boundary value problem by using fitted reproducing kernel algorithm. The reproducing kernel function is built to get fast accurately and efficiently series solutions with easily computable coefficients throughout evolution the algorithm under constraint periodic conditions within required grid points. The analytic solution is formulated in a finite series form whilst the truncated series solution is given to converge uniformly to analytic solution. The reproducing kernel procedure is based upon generating orthonormal basis system over a compact dense interval in Sobolev space to construct a suitable analytical-numerical solution. Furthermore, experiments result of some numerical examples are presented to illustrate the good performance of the presented algorithm. The results indicate that the reproducing kernel procedure is powerful tool for solving other problems of ordinary and partial differential equations arising in physics, computer and engineering fields.

**Keywords:** Boundary Value Problem, Error Estimation and Error Bound, Reproducing Kernel Theory

## Introduction

Periodic boundary value problem (PBVP) is an active research of modern mathematics that can be found naturally in different branches of applied sciences, physics and engineering (Gregu, 1987; Minh, 1998; Ashyralyev *et al.*, 2009). It has many applications due to the fact that a lot practical problems in mechanics, electromagnetic, astronomy and electrostatics may be converted directly to such PBVP. However, it is difficult generally to get a closed form solution for PBVP in terms of elementary functions, especially, for nonlinear and nonhomogeneous cases. So, PBVP has attracted much attention and has been studied by several authors (Kong and Wang, 2001; Chu and Zhou, 2006; Liu *et al.*, 2007; Zehour *et al.*, 2008; Yu and Pei, 2010; Abu Arqub and Al-Smadi, 2014). The purpose of this analysis is to develop analytical-numerical method for handling third-order, two-point PBVP with given periodic conditions by an application of the reproducing kernel theory. More specifically, consider the general form of third-order BVP:

$$y'''(t) = F(t, y(t), y'(t), y''(t)), \quad 0 \leq t \leq 1 \quad (1)$$

with periodic conditions:

$$y(0) - y(1) = 0, \; y'(0) - y'(1) = 0, \; y''(0) - y''(1) = 0 \quad (2)$$

where, $y \in y \in W_2^4[0, 1]$ is unknown function to be determined, $F(t, v_1, v_2, v_3)$ is continuous in $W_2^1[0, 1]$ as $v_i = v_i(t) \in W_2^4[0, 1]$, $0 \leq t \leq 1$, $-\infty < v_i < 1$, $i = 1, 2, 3$, which is linear or nonlinear term depending on the problem discussed.

The numerical solvability of BVPs with periodic conditions of different order has been pursued in literature. To mention a few, the existence and multiplicity of positive solutions for PBVP in the forms $u'''(x) = \alpha(x) f(x, u(x))$, $u'''(x) = f(x, u(x)) + \rho u(x)$ and $u'''(x) = \rho u''(x) + f(x, u(x)) + \alpha(x)$ have been studied, respectively, by Liu *et al.* (2007), Chu and Zhou (2006) and Yu and Pei (2010). However, Al-Smadi *et al.* (2014) have developed an iterative method for handling system of first-order PBVPs based on the RKHM. While Hopkins and Kosmatov (2007) have provided the existence of at least one positive solutions of PBVP in the form $u'''(x) = f(x, u(x), u'(x), u''(x))$. Further, Taylors decomposition method for solving linear periodic equation $u'''(x) = \alpha_1(x) u''(x) + \alpha_2(x) u'(x) + \alpha_3(x) u(x) + \alpha_4(x)$ numerically is proposed by Ashyralyev *et al.* (2009). The reproducing kernel method has widely applications to construct a numeric-analytic solution of different types of IVPs and BVPs. For example, see the work in (Geng and Cui, 2012; Al-Smadi *et al.*, 2012; 2013; Komashynska and Al-Smadi, 2014; Abu-Gdairi and Al-Smadi, 2013; Abu Arqub *et al.*, 2012; 2013; 2015; Geng and Qian, 2013; Bushnaq *et al.*, 2016; Ahmad *et al.*, 2016) for more information's

---

Corresponding Author: M. Al-Smadi, Department of Applied Science, Ajloun College, Al-Balqa Applied University, Ajloun 26816, Jordan. Email: mhm.smadi@yahoo.com

about RKHS method and scientific applications. But on the other aspect as well, several numerical schemes have been applied to solve IVPs and BVPs. For example, we refer to the work in (Komashynska *et al*., 2016a; 2016b; 2016c; Momani *et al*., 2014; Abu-Gdairi *et al*., 2015; Al-Smadi *et al*., 2015; 2016; Al-Smadi, 2013; Moaddy *et al*., 2015; El-Ajou *et al*., 2015; Abuteen *et al*., 2016).

The structure of this article is organized as follows. In section 2, reproducing kernel spaces are described to compute its reproducing kernels functions in which every function satis.es the periodic conditions. In section 3 and 4, the analytical-numerical solutions of Equation 1 and 2 as well an iterative method for obtaining these solutions are presented with series formula in the space $W_2^4$ [0, 1]. The n-term numerical solution is obtained to converge uniformly to analytic solution. In section 5, some numerical examples are simulated to check the reasonableness of our theory and to demonstrate the high performance of the presented algorithm. Conclusions are summarized in the last section.

## Construction of Reproducing Kernel Functions

In this section, we present some basic results and remarks in the reproducing kernel theory and its applications.

*Definition 1*

Let $W_2^4$ [0, 1] = {$y$| $y$, $y'$, $y''$, $y'''$ are absolutely continuous on [0, 1] such that $y^{(4)} \in L^2$ [0, 1] and $y(0)-y(1) = 0$, $y'(0)-y'(1) = 0$, $y''(0)-y''(1) = 0$}. On the other hand, let $\langle y(t), z(t) \rangle_{W_2^4}$ be the inner product in the space $W_2^4$ [0, 1], which is defined by:

$$\langle y(t), z(t) \rangle_{W_2^4} = \sum_{i=0}^{3} y^{(i)}(0) z^{(i)}(0) + \int_0^1 y^{(4)}(v) z^{(4)}(v) dv \quad (3)$$

and the norm is $\|y\|_{W_2^4} = \sqrt{\langle y(t), y(t) \rangle_{W_2^4}}$, where $y, z \in W_2^4[0,1]$.

*Lemma 1*

The reproducing kernel function $K_s(t)$ of the Hilbert space $W_2^4$ [0, 1] can be given by:

$$K_s(t) = \begin{cases} \sum_{i=0}^{7} a_i(s) t^i, & t \leq s \\ \sum_{i=0}^{7} b_i(s) t^i, & t > s \end{cases} \quad (4)$$

*Proof*

Using several integrations by parts of $\int_0^1 y^{(4)}(v) K_s^{(4)}(v) dv$ to obtain that:

$$\langle y(t), K_s(t) \rangle_{W_2^4} = \sum_{i=0}^{3} y^{(i)}(0) \left[ \partial_t^i K_s(0) + (-1)^i \partial_t^{7-i} K_s(0) \right] + \sum_{i=0}^{3} (-1)^{i+1} y^{(i)}(1) \partial_t^{7-i} K_s(1) + \int_0^1 y(v) \partial_t^8 K_s(v) dv$$

If $K_s(t) \in W_2^4$ [0, 1], then $K_s(0) = K_s(1)$, $\partial_t K_s(0) = \partial_t K_s(1)$, $\partial_t^2 K_s(0) = \partial_t^2 K_s(1)$, as well as if $y \in W_2^4$ [0, 1], then $y^{(i)}(0) = y^{(i)}(1)$, $i = 0, 1, 2$. Therefore:

$$\langle y(t), K_s(t) \rangle_{W_2^4} = \sum_{i=0}^{3} y^{(i)}(0) \left[ \partial_t^i K_s(0) + (-1)^i \partial_t^{7-i} K_s(0) \right] + \sum_{i=0}^{3} (-1)^{i+1} y^{(i)}(1) \partial_t^{7-i} K_s(1) + \int_0^1 y(v) \partial_t^8 K_s(v) dv + c_1 (y(0) - y(1)) + c_2 (y'(0) - y'(1)) + c_3 (y''(0) - y''(1))$$

For each $s, t \in$ [0, 1], assume $K_s(t)$ satisfy the following:

$$\partial_t^4 K_s(1) = 0, \partial_t^3 K_s(0) - \partial_t^4 K_s(0) = 0,$$
$$K_s(0) + \partial_t^7 K_s(0) + c_1 = 0,$$
$$\partial_t^7 K_s(1) + c_1 = 0, \partial_t^6 K_s(1) - c_2 = 0,$$
$$\partial_t K_s(0) - \partial_t^6 K_s(0) + c_2 = 0,$$
$$\partial_t^2 K_s(0) + \partial_t^5 K_s(0) + c_3 = 0,$$
$$\partial_t^5 K_s(1) + c_3 = 0$$

Thus, we have $\langle y(t), K_s(t) \rangle_{W_2^4} = \int_0^1 y(v) \partial_t^8 K_s(v) dv$. Also, assume $K_s(t)$ satisfy that:

$$\partial_t^8 K_s(t) = \delta(s-t), \delta \; dirac-delta \; function \quad (5)$$

so, $\langle y(t), K_s(t) \rangle_{W_2^4} = y(s)$.

Next, we give the expression form of $K_s(t)$, to do this, the auxiliary formula of Equation 5 is $\lambda^8 = 0$ and their real solutions are $\lambda = 0$ with multiplicity 8. Hence, let the form of $K_s(t)$ be as defined in Equation 4. On the other hand of Equation 5, let $K_s(t)$ satisfy $\partial_t^8 K_s(t+0) = \partial_t^m K_s(t-0)$, $m = 0, 1,..., 6$. By integrating $\partial_t^8 K_s(t) = \delta(s-t)$ from $s-\varepsilon$ to $s+\varepsilon$ with respect to $t$ as well letting $\varepsilon \to 0$, we have the jump degree of $\partial_t^7 K_s(t)$ at $s = t$ such that $\partial_t^7 K_s(t+0) - \partial_t^7 K_s(t-0) = -1$.

Through the last computational results the unknown coefficients $a_i(s)$ and $b_i(s)$, $i = 0, 1,..., 7$ of $K_s(t)$ in Equation 4 can be obtained. However, the representation form of these coefficients using Maple 13 software package are provided by:

$$a_0(s) = 1,$$
$$a_1(s) = \frac{1}{\alpha_1} s\left(9244 - 6300s - 50820s^2 - 12705s^3 + 205611s^4 - 203035s^5 + 58005s^6\right),$$
$$a_2(s) = \frac{1}{\alpha_2} s\left(-25200 + 2041264s - 3024100s^2 - 756025s^3 + 2570499s^4 - 806113s^5 + 5s^6\right),$$
$$a_3(s) = \frac{1}{\alpha_3} s\left(-1829520 - 27216900s + 72804784s^2, 54887415s^3 + 11205561s^4 - 76873s^5 + 363s^6\right),$$
$$a_5(s) = \frac{1}{\alpha_4} s\left(-1829520 - 27216900s + 72804784s^2, 54887415s^3 + 11205561s^4 - 76873s^5 + 363s^6\right),$$
$$a_5(s) = \frac{1}{\alpha_5} s\left(49346640 - 138124504s + 74703740s^2 + 18675935s^3 - 5054705s^4 + 462685s^5 - 9791s^6\right),$$
$$a_6(s) = \frac{1}{\alpha_6} s\left(146169244 - 145159740s - 1537460s^2 - 384365s^3 + 1388055s^4 - 504739s^5 + 29005s^6\right),$$
$$a_7(s) = \frac{1}{\alpha_7} \left(-292354444 + 292345200s + 6300s^2 + 50820s^3 + 12705s^4 - 205611s^5 + 203035s^6 - 58005s^7\right),$$

$$b_0(s) = 1 - \frac{1}{5040} s^7,$$
$$b_1(s) = \frac{1}{\alpha_1} s\left(9244 - 6300s - 50820s^2 - 12705s^3 + 205611s^4 + \frac{36542311}{180} s^5 + 58005s^6\right),$$
$$b_2(s) = \frac{1}{\alpha_2} s\left(-25200 + 2041264s - 3024100s^2 - 756025s^3 - \frac{34351126}{15} s^4 - 806443s^5 + 5s^6\right),$$
$$b_3(s) = \frac{1}{\alpha_3} s\left(-1829520 - 27216900s + 72804784s^2 + 18201196s^3 + 11205561s^4 - 76873s^5 + 363s^6\right),$$
$$b_4(s) = \frac{1}{\alpha_4} s\left(-1829520 - 27216900s + 219549660s^2 - 54887415s^3 + 11205561s^4 - 76873s^5 + 363s^6\right),$$
$$b_5(s) = \frac{1}{\alpha_5} s\left(49346640 + 154229940s + 74703740s^2 + 18675935s^3 - 5054705s^4 + 462685s^5, 9791s^6\right),$$
$$b_6(s) = \frac{1}{\alpha_6} s\left(-146185200 - 145159740s - 1537460s^2 - 384365s^3 + 1388055s^4 - 504739s^5 + 29005s^6\right),$$
$$b_7(s) = \frac{1}{\alpha_7} s\left(292345200 + 6300s + 50820s^2 + 12705s^3 - 205611s^4 + 203035s^5 - 58005s^6\right)$$

where, $\alpha_1 = 292354444$, $\alpha_2 = 1169417776$, $\alpha_3 = 10524759984$, $\alpha_4 = 42099039936$, $\alpha_5 = 70165066560$, $\alpha_6 = 210495199680$ and $\alpha_7 = 1473466397760$.

*Definition 2*

Let $W_2^1[0, 1] = \{y|\ y$ is absolutely continuous on $[0, 1]$ and $y' \in L^2[0, 1]\}$ and let the inner product $\langle y(t), z(t) \rangle_{W_2^1}$ be written as:

$$\langle y(t), z(t) \rangle_{W_2^1} = y(0)z(0) + \int_0^1 y'(v) z'(v) dv$$

whereas the norm $\|\cdot\|_{W_2^1}$ is given such that $\|y\|_{W_2^1} \sqrt{\langle y(t), y(t) \rangle_{W_2^1}}$, $y, z \in W_2^1[0,1]$.

However, Geng and Cui (2007) show that the reproducing kernel function $Q_s(t)$ of $W_2^1[0, 1]$ can be given by:

$$Q_s(t) = \begin{cases} 1+t & t \leq s, \\ 1+s & t > s, \end{cases}$$

# Representation of Analytical-Numerical Solutions

First, as in (Geng *et al.*, 2013; 2014; 2015; Al-Smadi and Altawallbeh, 2013; Abu Arqub, 2015), we transform the problem into a differential operator. To do so, we de.ne an operator $T$ from the space $W_2^4[0, 1]$ into $W_2^1[0, 1]$ such that $Ty(t) = y'''(t)$. Therefore, Equation 1 and 2 can be converted equivalently into following form:

$$Ty(t) = F(t, y(t), y'(t), y''(t)), 0 \leq t \leq 1 \tag{7}$$

with periodic conditions:

$$y(0) - y(1) = 0, y'(0) - y'(1) = 0, y''(0) - y''(1) = 0 \quad (8)$$

where, $y \in W_2^4[0,1]$ and $F(t, v_1, v_2, v_3) \in W_2^1[0, 1]$ for $v_i = v_i(t) \in W_2^4[0, 1]$, $\infty < v_i < \infty$, $i = 1, 2, 3$ and $0 \le t \le 1$.

*Lemma 2*

The operator $T$ is linear bounded operator from $W_2^4[0, 1]$ into $W_2^1[0, 1]$.

*Proof*

We want to show that $\|Ty\|_{W_2^1}^2 \le M \|y\|_{W_2^4}^2$, where $M > 0$. From Definition 2, we get that:

$$\|Ty\|_{W_2^1}^2 = \langle Ty(t), Ty(t)\rangle_{W_2^1} = [(Ty)(0)]^2 + \int_0^1 \left[(Ty)'(t)\right]^2 dt$$

Since $y(s) = \langle y(t), K_s(t)\rangle_{W_2^4}$, $(Ty)(s) = \langle y(t), (TK_s)(t)\rangle_{W_2^4}$ and $(Ty)'(s) = \langle y(t), (TK_s)'(t)\rangle_{W_2^4}$, we have by using Schwarz's inequality that:

$$|(Ty)(s)| = \left|\langle y(s), (TK_s)(s)\rangle_{W_2^4}\right| \le \|TK_s\|_{W_2^4} \|y\|_{W_2^4} = M_1 \|y\|_{W_2^4}$$

where $M_1 > 0$. Also, $|(Ty)'(s)| = \left|\langle y(s), (TK_s)'(s)\rangle_{W_2^4}\right| \le \|(TK_s)'\|_{W_2^4} \|y\|_{W_2^4} = M_2 \|y\|_{W_2^4}$, where $M_2 > 0$. Thus,

$$\|Ty\|_{W_2^1}^2 = [(Ty)(0)]^2 + \int_0^1 \left[(Ty)'(t)\right]^2 dt \le (M_1^2 + M_2^2)\|y\|_{W_2^4}^2.$$

The linearity part is obvious.

Now, to construct an orthonormal basis system $\{\bar{\psi}_i(t)\}_{i=1}^\infty$ of the space $W_2^4[0, 1]$, we set firstly $\varphi_i(t) = K_{t_i}(t)$ and $\psi_i(t) = T^*\varphi_i(t)$, where $\{t_i\}_{i=1}^\infty$ is dense set on compact $[0, 1]$ and $T^*$ is the adjoint of $T$. Thus, we have that $\langle y(t), \psi_i(t)\rangle_{W_2^4} = \langle y(t), T^*\varphi_i(t)\rangle_{W_2^4} = \langle Ty(t), \varphi_i(t)\rangle_{W_2^1} = Ty(t_i)$ $i = 1, 2, \ldots$.

*Theorem 1*

Let $\{t_i\}_{i=1}^\infty$ be a dense subset on $[0, 1]$, then $\{\psi_i(t)\}_{i=1}^\infty$ is complete basis system of $W_2^4[0, 1]$.

*Proof*

It is easy to note that $\psi_i(t) = T^*\varphi_i(t) = \langle T^*\varphi_i(t), K_s(t)\rangle_{W_2^4}$ $= \langle \varphi_i(t), T_t K_s(t)\rangle_{W_2^1} = T_t K_s(t)|_{y=x_i} \in W_2^4[0, 1]$. Therefore, $\psi_i(t)$ can be given by $\psi_i(t) = T_t K_s(t)|_{t=s_i}$, where $T_t$ refers to the operator $T$ applies to the function of $t$. However, let $\langle y(t), \psi_i(t)\rangle_{W_2^4} = 0$, $i = 1, 2, \ldots$ for each fixed $t \in W_2^4[0, 1]$. That is, $\langle y(t), \psi_i(t)\rangle_{W_2^4} = \langle y(t), T^*\varphi_i(t)\rangle_{W_2^4} = \langle Ty(t), \varphi_i(t)\rangle_{W_2^1} = Ty(t_i) = 0$. Since $\{t_i\}_{i=1}^\infty$ is dense on $[0, 1]$, thus $Ty(t) = 0$ which means that $y(t) = 0$.

Furthermore, $\{\bar{\psi}_i(t)\}_{i=1}^\infty$ can be obtained by using the Gram-Schmidt process of $\{\psi_i(t)\}_{i=1}^\infty$ as follows:

$$\bar{\psi}_i(t) = \sum_{k=1}^i \beta_{ik} \psi_k(t) \quad (9)$$

*Theorem 2*

Let $y(t)$ be analytic solution of Equation 7 and 8 and let $\{t_i\}_{i=1}^\infty$ be a dense subset on $[0, 1]$, then:

$$y(t) = T^{-1}F(t, y(t), y'(t), y''(t))$$
$$= \sum_{i=1}^\infty \sum_{k=1}^i \beta_{ik} F(t_k, y(t_k), y'(t_k), y''(t_k)) \bar{\psi}_i(t) \quad (10)$$

*Proof*

Clearly that $\sum_{i=1}^\infty \langle y(t), \bar{\psi}_i(t)\rangle \bar{\psi}_i(t)$ is Fourier series expansion about orthonormal basis $\{\bar{\psi}_i(t)\}_{i=1}^\infty$, $y \in W_2^4[0, 1]$. Anyhow, since $W_2^4[0, 1]$ is the Hilbert space, then $\sum_{i=1}^\infty \langle y(t), \bar{\psi}_i(t)\rangle \bar{\psi}_i(t)$ is convergent in the sense of $\|\cdot\|_{W_2^4}$. Whilst by Equation 9, we have that:

$$y(t) = \sum_{i=1}^\infty \langle y(t), \bar{\psi}_i(t)\rangle_{W_2^4} \bar{\psi}_i(t)$$
$$= \sum_{i=1}^\infty \sum_{k=1}^i \beta_{ik} \langle y(t), \psi_k(t)\rangle_{W_2^4} \bar{\psi}_i(t)$$
$$= \sum_{i=1}^\infty \sum_{k=1}^i \beta_{ik} \langle y(t), T^*\varphi_k(t)\rangle_{W_2^4} \bar{\psi}_i(t)$$
$$= \sum_{i=1}^\infty \sum_{k=1}^i \beta_{ik} \langle Ty(t), \varphi_k(t)\rangle_{W_2^1} \bar{\psi}_i(t)$$
$$= \sum_{i=1}^\infty \sum_{k=1}^i \beta_{ik} \langle F(t, y(t), y'(t), y''(t)), \varphi_k(t)\rangle_{W_2^1} \bar{\psi}_i(t)$$
$$= \sum_{i=1}^\infty \sum_{k=1}^i \beta_{ik} F(t_k, y(t_k), y'(t_k), y''(t_k)) \bar{\psi}_i(t)$$
$$= T^{-1}F(t, y(t), y'(t), y''(t))$$

Let $\{\bar{\psi}_i(t)\}_{i=1}^{\infty}$ be the orthonormal basis obtained by Gram-Schmidt process of $\{\psi_i(t)\}_{i=1}^{\infty}$, then the analytic solution of Equation 7 and 8 can be written as follows:

$$y(t) = \sum_{i=1}^{\infty} B_i \bar{\psi}_i(t) \qquad (11)$$

where, $\sum_{k=1}^{i} \beta_{ik} F(t_k, y(t_k), y'(t_k), y''(t_k))$. Indeed, $B_i$ are unknown. Thus, we can approximate $B_i$ using known $A_i$. For computations, define the initial guess function $y_0(t_1) = 0$, set $y_0(t_1) = y(t_1)$ and define the $n$th-order approximation $y_n(t)$ to $y(t)$ as follows:

$$y_n(t) = \sum_{i=1}^{n} A_i \bar{\psi}_i(t) \qquad (12)$$

where, the coefficients $A_i$ of $\bar{\psi}_i(t)$, $i = 1, 2, ..., n$ are obtained by:

$$A_i = \sum_{k=1}^{i} \beta_{ik} F(t_k, y_{k-1}(t_k), y'_{k-1}(t_k), y''_{k-1}(t_k)) \qquad (13)$$

*Theorem 3*

If $y \in W_2^4[0,1]$, then there exists $M > 0$ such that $\|y^{(i)}\|_C \leq M \|y\|_{W_2^4}$, $i = 0, 1, 2, 3$, where $\|y\|_C = \max_{0 \leq t \leq 1} |y(t)|$.

*Proof*

For each $s, t \in [0, 1]$, we have $y^{(i)}(t) = \langle y(t), \partial_t^i K_s(t) \rangle_{W_2^4}$, $i = 0, 1, 2, 3$. By the expression form of $K_s(t)$, it follows that $\|\partial_t^i K_s\|_{W_2^4} \leq M_i$, $i = 0, 1, 2, 3$, where $M_i > 0$. Thus, $|y^{(i)}(t)| = |\langle y(t), \partial_t^i K_s(t) \rangle_{W_2^4}| \leq \|\partial_t^i K_s\|_{W_2^4} \|y\|_{W_2^4} \leq M_i \|y\|_{W_2^4}$, $i = 0, 1, 2, 3$. Hence, $\|y^{(i)}\|_C \leq \max_{i=0,1,2,3} \{M_i\} \|y\|_{W_2^4}$, $i = 0, 1, 2, 3$.

*Corollary 1*

The numeric solution and its derivatives up to order three are converge uniformly to analytic solution and all its derivatives, respectively.

*Proof*

For each $t \in [0, 1]$, we have:

$$|y_n^{(i)}(t) - y^{(i)}(t)| = |\langle y_n(t) - y(t), \partial_t^i K_s(t) \rangle_{W_2^4}|$$
$$\leq \|\partial_t^i K_x\|_{W_2^4} \|y_n - y\|_{W_2^4}$$
$$\leq M_i \|y_n - y\|_{W_2^4}, i = 0, 1, 2, 3$$

where, $M_i > 0$. Therefore, if $\|y_n - y\|_{W_2^4} \to 0$ as $n \to 1$, then the numeric solution $y_n(x)$ and $y_n^{(i)}(t)$, $i = 1, 2, 3$, are converge uniformly to analytic solution $y(t)$ and $y^{(i)}(t)$, $i = 1, 2, 3$, respectively.

## Convergence Analysis and Error Estimation

From Corollary 1, if $y_n(t)$ converging uniformly to $y(t) = \sum_{i=1}^{\infty} A_i \bar{\psi}_i(t)$. If $y_n(t) = P_n y(t)$, where $P_n$ is orthogonal project from the space $W_2^4[0, 1]$ to Span $\{\psi_1, \psi_2, ..., \psi_n\}$, then $T y_n(t_j) = \langle T y_n(t), \varphi_j(t) \rangle_{W_2^1} = \langle y_n(t), T_j^* \varphi(t) \rangle_{W_2^4} = \langle P_n y(t), \psi_j(t) \rangle_{W_2^4} = \langle y(t), P_n \psi_j(t) \rangle_{W_2^4} = \langle y(t), \psi_j(t) \rangle_{W_2^4} = \langle Ty(t), \varphi_j(t) \rangle_{W_2^1} = Ty(t_j)$. Next, we list two lemmas for convenience in order to prove the recent theorems.

*Lemma 3*

The numerical solution $y_n$ satisfies, $T y_n(t_j) = F(t_j, y_{j-1}(t_j), y'_{j-1}(t_j), y''_{j-1}(t_j))$ for $j = 1, 2, 3,$.

*Proof*

Let $A_i = \sum_{k=1}^{i} \beta_{ik} F(t_k, y_{k-1}(t_k), y'_{k-1}(t_k), y''_{k-1}(t_k))$, then $y_n$ will be rewritten in form of $y_n(t) = \sum_{i=1}^{n} A_i \psi_i(t)$. Using the properties of $K_s(t)$, it follows that:

$$Ty_n(t_j) = \sum_{i=1}^{n} A_i T \psi_i(t_j) = \sum_{i=1}^{n} A_i \langle T \bar{\psi}_i(t) \rangle_{W_2^1}$$
$$= \sum_{i=1}^{n} A_i \langle \bar{\psi}_i(t), T_j^* \varphi(t) \rangle_{W_2^4} = \sum_{i=1}^{n} A_i \langle \bar{\psi}_i(t), \psi_j(t) \rangle_{W_2^4}$$

By orthogonality of $\{\bar{\psi}_i(t)\}_{i=1}^{\infty}$, we have:

$$\sum_{l=1}^{j} \beta_{jl} T y_n(t_l) = \sum_{i=1}^{n} A_i \left\langle \bar{\psi}_i(t), \sum_{l=1}^{j} \beta_{jl} \psi_l(t) \right\rangle_{W_2^4}$$
$$= \sum_{i=1}^{n} A_i \langle \bar{\psi}_i(t), \bar{\psi}_j(t) \rangle_{W_2^4}$$
$$A_j = \sum_{l=1}^{j} \beta_{jl} F(t_l, y_{l-1}(t_l), y'_{l-1}(t_l), y''_{l-1}(t_l))$$

Taking $j = 1$, one gets $T y_n(t_1) = F(t_1, y_0(t_1), y'_1(t_1), y''_0(t_1))$. Taking $j = 2$, one gets:

$$Ty_n(t_2) = F(t_2, y_1(t_2), y'_1(t_2), y''_1(t_2))$$

Therefore via mathematical induction, we can get that $Ty_n(t_j) = F(t_j, y_{j-1}(t_j), y'_{j-1}(t_j), y''_{j-1}(t_j))$.

*Lemma 4*

Let $F(t, v_1, v_2; v_3)$ be continuous function in [0, 1] with respect to $t$, $v_i$, $t \in [0, 1]$, where $v_i \in (-\infty, \infty)$ for $i = 1, 2, 3$. If $\|y_n - y\|_{W_2^4} \to 0$, $t_n \to t$ as $n \to 1$, then $F(t_n, y_{n-1}(t_n), y'_{n-1}(t_n), y''_{n-1}(t_n)) \to F(t, y(t), y'(t), y''(t))$ as $n \to 1$.

*Proof*

Since $\|y_n - y\|_{W_2^4} \to 0$ as $n \to 1$, by Corollary .1, it follows that $y_{n-1}^{(i)}(t)$ is converging uniformly to $y^{(i)}(t)$, $i = 0, 1, 2$. Hence, the proof is completed since F is continuous.

*Theorem 4*

Let $\{t_i\}_{i=1}^\infty$ be dense on [0, 1] and $\|y_n\|_{W_2^3}$ be bounded, then $y_n(t)$ in Equation 12 converges to $y(t)$ of Equation 7 and 8 in the space $W_2^4[0, 1]$ such that $y(t) = \sum_{i=1}^\infty A_i \bar{\psi}_i(t)$, where $A_i$ is obtained in Equation 13.

*Proof*

Firstly, we will show that $\{y_n\}_{i=1}^\infty$ in Equation 12 is increasing by sense of the norm of $W_2^4[0, 1]$. Since $\{\bar{\psi}_i\}_{i=1}^\infty$ is complete normal orthogonal basis in $W_2^4[0, 1]$, then

$$\|y_n\|_{W_2^4}^2 = \langle y_n(t), y_n(t) \rangle_{W_2^4} = \left\langle \sum_{i=1}^n A_i \bar{\psi}_i(t), \sum_{i=1}^n A_i \bar{\psi}_i(t) \right\rangle_{W_2^4} = \sum_{i=1}^n (A_i)^2$$

. Therefore, $\|y_n\|_{W_2^4}^2$ is increasing.

Secondly, we will show a convergence of $y_n(t)$. By Equation (12), $y_{n+1}(t) = y_n(t) + A_{n+1}\bar{\psi}_{n+1}(t)$. Hence, we have $\|y_{n+1}\|_{W_2^4}^2 = \|y_n\|_{W_2^4}^2 + (A_{n+1})^2 = \|y_{n-1}\|_{W_2^4}^2 + (A_n)^2 + (A_{n+1})^2 = \ldots = \|y_0\|_{W_2^4}^2 + \sum_{i=1}^{n+1} (A_i)^2$. Since, the sequence $\|y_n\|_{W_2^4}^2$ is increasing as well as $\|y_n\|_{W_2^4}$ is bounded, then $\|y_n\|_{W_2^4}$ is convergent as $n \to 1$. That is, $\exists$ a constant $\alpha$ such that $\sum_{i=1}^\infty (A_i)^2 = \alpha$. Thus $A_i = \sum_{k=1}^j \beta_{ik} F(t_k, y_{k-1}(t_k), y'_{k-1}(t_k), y''_{k-1}(t_k)) \in l^2$, $i = 1, 2, \ldots$. On the other hand, since $(y_m - y_{m-1}) \perp (y_{m-1} - y_{m-2}) \perp \ldots \perp (y_{n+1} - y_n)$ it implies that $\|y_m - y_n\|_{W_2^4}^2 = \|y_m - y_{m-1} + y_{m-1} - \ldots + y_{n+1} - y_n\|_{W_2^4}^2$
$= \|y_m - y_{m-1}\|_{W_2^4}^2 + \ldots + \|y_{n+1} - y_n\|_{W_2^4}^2$ for $m > n$. Furthermore, $\|y_m - y_{m-1}\|_{W_2^4}^2 = (A_m)^2$. Anyhow, $\|y_m - y_n\|_{W_2^4}^2 = \sum_{i=n+1}^m (A_i)^2 \to 0$ as $n, m \to 1$. From completeness of $W_2^4[0, 1]$, for $n \to 1$, $\exists y \in W_2^4[0, 1]$ such that $y_n(t) \to y(t)$ in the sense of $\|\cdot\|_{W_2^4}$. Finally, we will show that $y(t)$ is analytic solution of Equation 7 and 8.

Let $\{t_i\}_{j=1}^\infty$ be dense on [0, 1], then there is a subsequence $\{t_{n_j}\}_{j=1}^\infty$ such that $t_{n_j} \to t$ as $j \to \infty$, $\forall t \in [0, 1]$. It is clear, by Lemma .3, that $Ty(t_{n_j}) = F(t_{n_j}, y_{n_j-1}(t_k), y'_{n_j-1}(t_k), y''_{n_j-1}(t_k))$. From lemma .4 and continuity of F, it follows for $j \to \infty$ that $Ty(t) = F(t, y(t); y'(t), y''(t))$. That is, $y(t)$ satisfies Equation 7. As well as if $\bar{\psi}_i(t) \in W_2^4[0, 1]$, then $y(t)$ satisfy periodic conditions of Equation 8, which means that $y(t)$ is analytic solution of Equation 7 and 8 such that $y(t) = \sum_{i=1}^\infty A_i \bar{\psi}_i(t)$.

# Computational Algorithm and Numerical Experiments

Using RKHS method, taking $t_i = \frac{i-1}{n-1}$, $i = 1, 2, \ldots, n$ and according to reproducing kernel functions $K_s(t)$ and $Q_s(t)$ on [0, 1]; some tabulate results are presented quantitatively at some selected grid points on [0, 1] to illustrate the accuracy of the FRK method for handling the PBVP.

*Example 1*

Meditate in the following linear differential equation:

$$y'''(t) + ty''(t) f(t), \ 0 \le t \le 1$$

with periodic conditions:

$$y(0) - y(1) = 0, y'(0) - y'(1) = 0, y''(0) - y''(1) = 0$$

where, $f(t)$ is given to obtain the exact solution as $y(t) = e^{t^2(1-t)^2}$.

*Example 2*

Meditate in the following nonlinear diferential equation:

$$y'''(t) - \cos(t)(y''(t))^3 + 2\sinh(y(t))\cosh(y(t)y'(t))$$
$$= f(t), \ 0 \le t \le 1$$

with periodic conditions:

$$y(0) - y(1) = 0, y'(0) - y'(1) = 0, y''(0) - y''(1) = 0$$

where, $f(t)$ is given to obtain the exact solution as $y(t) = \frac{1}{2}t^2(1-t)^2$.

*Example 3*

Meditate in the following nonlinear differential equation:

$$y'''(t) + y''(t) + t(y'(t))^2 - \cosh^{-1}(y(t)) = f(t), 0 \leq t \leq 1$$

with periodic conditions:

$$y(0) - y(1) = 0, y'(0) - y'(1) = 0, y''(0) - y''(1) = 0$$

where, $f(t)$ is given to obtain the exact solution as $y(t) = \cosh(t^2-t)$.

The agreement between the analytical-numerical solutions is investigated for Examples 1, 2 and 3 at various $t$ in [0, 1] by computing absolute errors and relative errors of numerically approximating their analytical solutions as shown in Tables 1 to 3, respectively.

Table 1. The analytical-numerical solutions and errors for Example 1

| $t$ | $y(t)$ | $y_{51}(t)$ | $|y(t)-y_{51}(t)|$ | $|y(t)|^{-1} |y(t)-y_{51}(t)|$ |
|---|---|---|---|---|
| 0 | 1 | 1 | 0 | 0 |
| 0.1 | 1.0081328937531524 | 1.0081340879367722 | $1.19418\times10^{-6}$ | $1.18455\times10^{-6}$ |
| 0.2 | 1.0259304941903822 | 1.0259325793058158 | $2.08512\times10^{-6}$ | $2.03241\times10^{-6}$ |
| 0.3 | 1.0450868583490185 | 1.0450883361164949 | $1.47777\times10^{-6}$ | $1.41401\times10^{-6}$ |
| 0.4 | 1.0592911944779007 | 1.0592914256899810 | $2.31212\times10^{-7}$ | $2.18271\times10^{-7}$ |
| 0.5 | 1.0644944589178593 | 1.0644940798919150 | $3.79026\times10^{-7}$ | $3.56062\times10^{-7}$ |
| 0.6 | 1.0592911944779007 | 1.0592914066684471 | $2.12191\times10^{-7}$ | $2.00314\times10^{-7}$ |
| 0.7 | 1.0450868583490185 | 1.0450883036692906 | $1.44532\times10^{-6}$ | $1.38297\times10^{-6}$ |
| 0.8 | 1.0259304941903820 | 1.0259325450574925 | $2.05087\times10^{-6}$ | $1.99903\times10^{-6}$ |
| 0.9 | 1.0081328937531522 | 1.0081340662518540 | $1.17250\times10^{-6}$ | $1.16304\times10^{-6}$ |
| 1 | 1 | 1 | 0 | 0 |

Table 2. The analytical-numerical solutions and errors for Example 2

| $t$ | $y(t)$ | $y_{36}(t)$ | $|y(t)-y_{36}(t)|$ | $|y(t)|^{-1} |y(t)-y_{36}(t)|$ |
|---|---|---|---|---|
| 0 | 0 | 0 | 0 | Indeterminate |
| 0.1 | 0.00405 | 0.004049976999474401 | $2.30005\times10^{-8}$ | $5.67914\times10^{-6}$ |
| 0.2 | 0.01280 | 0.012799927364010002 | $7.26360\times10^{-8}$ | $5.67469\times10^{-6}$ |
| 0.3 | 0.02205 | 0.022049874903704060 | $1.25096\times10^{-7}$ | $5.67330\times10^{-6}$ |
| 0.4 | 0.02880 | 0.028799836624554605 | $1.63375\times10^{-7}$ | $5.67276\times10^{-6}$ |
| 0.5 | 0.03125 | 0.031249822731293803 | $1.77269\times10^{-7}$ | $5.67260\times10^{-6}$ |
| 0.6 | 0.02880 | 0.028799836624554600 | $1.63375\times10^{-7}$ | $5.67276\times10^{-6}$ |
| 0.7 | 0.02205 | 0.022049874903704060 | $1.25096\times10^{-7}$ | $5.67330\times10^{-6}$ |
| 0.8 | 0.01280 | 0.012799927364010000 | $7.26360\times10^{-8}$ | $5.67469\times10^{-6}$ |
| 0.9 | 0.00405 | 0.004049976999474354 | $2.30005\times10^{-8}$ | $5.67914\times10^{-6}$ |
| 1 | 0 | 0 | 0 | Indeterminate |

Table 3. The analytical-numerical solutions and errors for Example 3

| $t$ | $y(t)$ | $y_{26}(t)$ | $|y(t)-y_{26}(t)|$ | $|y(t)|^{-1} |y(t)-y_{26}(t)|$ |
|---|---|---|---|---|
| 0 | 1 | 1 | 0 | 0 |
| 0.1 | 1.0040527344882193 | 1.0040527445966834 | $1.01085\times10^{-8}$ | $1.00677\times10^{-8}$ |
| 0.2 | 1.0128273299790107 | 1.0128269444788909 | $3.85500\times10^{-7}$ | $3.80618\times10^{-7}$ |
| 0.3 | 1.0221311529634651 | 1.0221299798303292 | $1.17313\times10^{-6}$ | $1.14773\times10^{-6}$ |
| 0.4 | 1.0289385056939790 | 1.0289365747582975 | $1.93094\times10^{-6}$ | $1.87663\times10^{-6}$ |
| 0.5 | 1.0314130998795732 | 1.0314108598211662 | $2.24006\times10^{-6}$ | $2.17183\times10^{-6}$ |
| 0.6 | 1.0289385056939790 | 1.0289365727521844 | $1.93294\times10^{-6}$ | $1.87858\times10^{-6}$ |
| 0.7 | 1.0221311529634651 | 1.0221299763027610 | $1.17666\times10^{-6}$ | $1.15118\times10^{-6}$ |
| 0.8 | 1.0128273299790107 | 1.0128269404245114 | $3.89554\times10^{-7}$ | $3.84621\times10^{-7}$ |
| 0.9 | 1.0040527344882193 | 1.0040527415434230 | $7.05520\times10^{-9}$ | $7.02673\times10^{-9}$ |
| 1 | 1 | 1 | 0 | 0 |

## Concluding Summary

In this article, we introduce the fitting reproducing kernel approach to enlarge its application range for treating a class of third-order periodic BVPs in a favorable reproducing kernel Hilbert space. The method does not require discretization of the variables as well as it provides best solutions in a less number of iterations and reduces the computational work. Further, we can conclude that the presented method is powerful and efficient technique in finding approximate solution for both linear and nonlinear problems. In the proposed algorithm, the solution and its approximation are represented in the form of series in $W_2^4[0, 1]$. The approximate solution and its derivative converge uniformly to exact solution and its derivative, respectively.

## Acknowledgment

The authors express their thanks to unknown referees for the careful reading and helpful comments.

## Author's Contributions

All authors completed the paper together as well as read and approved the final manuscript.

## Ethics

The authors declare that there is no conflict of interests regarding the publication of this paper.